\documentclass[12pt]{article}

\usepackage{amsmath, amsthm, amsfonts, amssymb}
\usepackage{booktabs,longtable,array}
\usepackage{caption,graphicx}
\newcommand{\tern}[3]{\langle #1,#2,#3\rangle}
\makeatletter
\newcommand{\namedlabel}[2]{\protected@edef\@currentlabel{#2}\label{#1}}
\makeatother

\newtheorem{prob}{Problem}
\newtheorem{lemma}{Lemma}[section]
\newtheorem{theorem}[lemma]{Theorem}
\newtheorem{proposition}[lemma]{Proposition}
\theoremstyle{definition}
\newtheorem{definition}[lemma]{Definition}
\newtheorem{example}[lemma]{Example}
\newtheorem{remark}[lemma]{Remark}

\newcommand{\N}{\mathbb{N}}
\newcommand{\Z}{\mathbb{Z}}
\newcommand{\lp}{\longrightarrow}

\title{Universality of monomial square forms}
\author{Stefan Maubach\\ \small Maastricht University, Maastricht, The Netherlands}
\date{}

\begin{document}

\maketitle

\begin{abstract} A monomial square form is a polynomial $F=\sum_{i=1}^k M_i^2$, where $M_1,\ldots,M_k$ are pairwise distinct squarefree monomials $\ne1$. 
We consider the problem of which of these forms are universal, i.e.\ represent every positive integer. 
We consider three classes of such forms:\\
(1) Forms $M_1^2+M_2^2+M_3^2+M_4^2$ ($k=4$),\\
(2) Forms $M_1^2+\ldots +M_5^2$ ($k=5$) in at most 5 variables,\\
(3) Forms where the monomials are products of consecutive variables of a fixed length:\\
$P_{m,d}(x)=\sum_{i=1}^{d}(x_ix_{i+1}\cdots x_{i+m-1})^2,$ and its circular version\\
$C_{m,d}(x)=\sum_{i=1}^{n}(x_ix_{i+1}\cdots x_{i+m-1})^2\quad(n=m+d-1,\ x_{i+n}=x_i).$
We classify which of these are universal, possibly universal, or conditionally universal on a conjecture on the 
form $x^2+4y^2+5z^2$. This ternary form appears repeatedly and seems to play a central role. 

For forms (1), we also show that a form is universal if and only if it represents 14, 15, 29 and 9591 (partially conditional on $x^2+4y^2+5z^2$). 
\end{abstract}



\tableofcontents
\section{Preliminaries}

\subsection{Introduction}
In this article we consider monomials $M_i$ of the form
\[ M_i= x_1^{e_1}x_2^{e_2}\cdots x_n^{e_n} , ~~ e_i\in \{0,1\}.\]
We disregard the case where all $e_i=0$ (the monomial 1) as that case is easily handled with the other cases in this article. 
  We consider forms  
\[ F=\sum_{i=1}^{k} M_i^2, ~M_i\textup{~pairwise~distinct}  \]
i.e.\ sums of squares of such monomials, each with coefficient~$1$; in particular $F(0,\ldots,0)=0$. The requirement that the $M_i$ be pairwise distinct ensures that all coefficients are~$1$ (two equal monomials would combine to $2M^2$).

Equivalently, these are the polynomials with zero constant term and all coefficients equal to~$1$ in which every variable occurs only to the powers $0$ and $2$. We call such forms {\bf monomial square forms}. 

In this article, we will only consider integer inputs for variables of forms.
On integer inputs such a form $\sum_{i=1}^{k} M_i^2$ can yield only nonnegative integers. If the form represents every positive integer, the form is called {\bf universal}.  This is the connection with quadratic forms: homogeneous polynomials of degree 2. The integers represented by quadratic forms have been studied extensively~\cite{Cassels,Dickson}; this is a central topic in number theory.

Monomial square forms specialize to diagonal quadratic forms (plus potentially a constant) by substituting constants for some of the variables. From this viewpoint, the set of values of a monomial square form is a union of value sets of (shifted) diagonal quadratic forms. 

In this article, we will be mainly concerned with when a monomial square form is universal. When analyzing this, there are two main parameters that order forms: the number of monomials, and within that the number of variables. In this article, we will consider three cases:
\begin{itemize}
\item At most four monomials, in any number of variables (Section~\ref{4less}).
\item Five monomials in at most five variables (Section~\ref{5m5v}).
\item Two special families of forms, in any number of variables. (Section~\ref{PnmCnm})
The forms are 
\begin{gather*} P_{m,d}(x)=\sum_{i=1}^{d}(x_ix_{i+1}\cdots x_{i+m-1})^2,\\
C_{m,d}(x)=\sum_{i=1}^{n}(x_ix_{i+1}\cdots x_{i+m-1})^2\quad(x_{i+n}=x_i),\end{gather*}
where $m\ge2$, $d\ge2$ and $n=m+d-1$: each monomial has $m$ variables, $P_{m,d}$ has $d$ monomials, and $C_{m,d}$ has $n$.
\end{itemize}

\subsection{Motivation}

First of all, monomial square forms are in some sense a very natural extension of quadratic forms, sharing certain properties of them (in particular being nonnegative and of degree $0$ or $2$ in each variable). Here we only consider monomial square forms with coefficients 1. 

Another important motivation for this article is pure curiosity. A second aim is to find challenging new problems that seem to need new ideas and do not yield to a routine (AI-assisted) computation. This is why we included Section~\ref{s:open}, where we state the most interesting open problems so that they can be read independently of the rest of the paper. 

A motivation for considering specifically the $P_{m,d}$ and $C_{m,d}$ forms is that 
they are in some sense further away from quadratic forms than many other forms considered in this article, as they ``suffer'' the most from sharing variables among their monomials. 
Also, they are forms with many symmetries (especially the $C_{m,d}$). Studying them is also an example of how one could make use of the results of Sections~\ref{4less} and~\ref{5m5v}. 

\subsection{Acknowledgements}

We would like to thank Prof. Andrew Earnest for some very useful email exchanges. 

There has been extensive use of Claude AI (mainly Opus 5.5) in the preparation of this article. Perhaps nowadays one should explain what the human contribution is, instead of explaining what the AI contribution is: the calculations are done by AI but the leading of the project and picking what is important has been a human endeavor. Also, many key ideas are human (like some of the lemmas) though the application of them has been mostly by AI in its computations. The fact that this article is somewhat readable is definitely due to a human (the author). 

We must say that the fact that such an AI model exists, has caused the author to do research again, as he's a lazy bugger and likes the fact that he doesn't need to do as much elaborate calculations etc. 

\subsection{Conventions}\label{conv}

(Unfortunately) we need to introduce quite a number of definitions and terminology, but most of them are quite natural.
We write $\N=\{0,1,2,\ldots\}$ and $\N^*=\{1,2,3,\ldots\}$.
Since only squares of variables occur, $F(\Z^n)=F(\N^n)$, and we work with $F(\N^n)$ throughout.

\begin{definition}We introduce the following terms for polynomials in $\Z[x_1,\ldots,x_n]$:
\begin{itemize}
\item A {\bf renaming} of a polynomial $P$ is constructed by renaming the variables appearing in $P$. 
A {\bf renaming class} is a set of polynomials that are renamings of each other. In Sections~\ref{4less} and \ref{5m5v} and the appendices, ``class'' always means renaming class.
\item A {\bf specialization} of a polynomial $P$ is a polynomial $P'$ obtained by putting some variables appearing in $P$ equal to integer constants.
\item A {\bf renamed specialization} of a polynomial $P$ is a specialization $P'$ of a renaming of $P$; we write $P\lp P'$ and also say that $P$ \emph{specializes to} $P'$. 
\item A {\bf value-preserving specialization} of a polynomial $P$ is a polynomial $P'$ such that $P\lp P'$ and $P$ and $P'$ have the same set of values. (In particular, every renaming is a value-preserving specialization.)
We write $P\twoheadrightarrow P'$. 
\item A polynomial $P$ is {\bf minimal} if every value-preserving specialization of $P$ is a renaming.
Minimal forms are the ``essential'' objects: they cannot be simplified without losing values.
For example, the forms $x^2+(yz)^2$ and $(xyz)^2+(yzt)^2$ both have $x^2+y^2$ as a value-preserving specialization (put $z=1$, respectively $y=z=1$ and rename), and $x^2+y^2$ is minimal. See also Lemma~\ref{l:types} below.

If $P\lp F$, $F$ is universal and $F(0,\ldots,0)=0$, then $P$ is universal too and $P\twoheadrightarrow F$.

\item Let $F=\sum_{i=1}^k M_i^2$ be a monomial square form. The {\bf fingerprint} of a
variable $x$ of $F$ is the vector $f(x)=(f_1(x),\ldots,f_k(x))\in\{0,1\}^k$, where
$f_i(x)=1$ if $x$ divides $M_i$ and $f_i(x)=0$ otherwise.
 The fingerprints determine $F$, since
$M_i=\prod_{f_i(x)=1}x$; renumbering the monomials permutes the coordinates.
If two variables $x,y$ have the same fingerprint, they only occur in $F$ through the
product $xy$, so putting $y=1$ is a value-preserving specialization. Hence every minimal
form has variables with pairwise distinct fingerprints, and every form has a
value-preserving specialization to such a form (Lemma~\ref{l:types}). Not every such
form is minimal: the variables of $(xu)^2+y^2+z^2+(tu)^2$ have pairwise distinct
fingerprints, but $u=1$ gives $x^2+y^2+z^2+t^2$, and both forms are universal.

\end{itemize}
\end{definition}

\begin{definition} For a form $F$, we define $E(F)$ to be the set of positive integers not represented by $F$.
\end{definition}

\begin{definition} 
We will use the following notation to represent a quadratic form $a_1x_1^2+\ldots+a_nx_n^2$:
\[ \langle a_1,\ldots,a_n\rangle. \]
Renaming is the same as putting the $a_i$ in a different order, and we can assume $a_1\leq a_2\leq \ldots \leq a_n$ up to renaming.

We will also write $c+a_1x_1^2+\ldots+a_nx_n^2$ as $c+\langle a_1,\ldots,a_n\rangle$ where $c$ is  a constant.
\end{definition}

\begin{definition}\label{def:branch}
Let $F$ be a monomial square form. If we have a specialization of $F$ which becomes 
$ c+\langle a_1,\ldots,a_n\rangle$, then we call this a {\bf branch} of $F$. 
We will call $c$ its {\bf shift}. 
\end{definition}
Since a specialization can only lose values, if $c+Q$ is a branch of $F$ and $n-c$ is represented by $Q$, then $n$ is represented by $F$. Very often the universality of $F$ is decided by finding appropriate ternary (or sometimes quaternary) branches; a branch is \emph{ternary} if its quadratic form has three variables.

\section{Tools used}

\label{Tools}

In this section we list known results and arguments we use in the following sections. 

\subsection{Known theorems used and a conjecture}

\begin{itemize}
\item[{\bf T15}]\namedlabel{15}{T15}(The 15-theorem) (Conway--Schneeberger 1993~\cite{Conway}; proof by Bhargava~\cite{Bhargava}). A positive definite quadratic form with integer matrix, i.e.\ $Q(x)=\sum_{i,j}a_{ij}x_ix_j$ with $a_{ij}=a_{ji}\in\mathbb Z$, represents every positive integer if and only if it represents each of
\[ 1,\ 2,\ 3,\ 5,\ 6,\ 7,\ 10,\ 14,\ 15. \]
\item[{\bf (FE)}]\namedlabel{FE}{(FE)} (Fermat--Euler~\cite{HW}) $n\ge1$ is a sum of two squares if and only if every prime $q\equiv3\pmod4$ divides $n$ to an even power.
\item[{\bf T1111}]\namedlabel{T1111}{T1111} (Lagrange~\cite{HW}) The form $\langle 1,1,1,1\rangle$ is universal.
\item[{\bf T111}]\namedlabel{T111}{T111} (Gauss--Legendre~\cite{Gauss}) $\langle 1,1,1 \rangle $ represents exactly the $n\ge0$ not of the form $4^a(8b+7)$.
\item[{\bf T112}]\namedlabel{T112}{T112} (Dickson~\cite{Dickson1927,JP}) $\langle 1,1,2\rangle $ represents exactly the $n\ge0$ not of the form $4^a(16b+14)$.
\item[\bf T122]\namedlabel{T122}{T122} (Dickson, Jones--Pall~\cite{Dickson1927,JP}) $\langle 1,2,2\rangle$ represents exactly the $n\ge0$ not of the form $4^a(8b+7)$.
\item[\bf T114]\namedlabel{T114}{T114} (Lemma~\ref{l:d}) $\langle 1,1,4\rangle$ represents every positive $n\not\equiv3\pmod4$ not of the form $4^a(8b+7)$.
\item[\bf T144]\namedlabel{T144}{T144} (Lemma~\ref{l:d}) $\langle 1,4,4\rangle$ represents every positive $n\equiv1\pmod4$, and every $n=4n'$ with $n'$ not of the form $4^a(8b+7)$.
\item[\bf T11616]\namedlabel{T11616}{T11616} (Jones--Pall~\cite{JP}) Write $n=4^au$ with $4\nmid u$. Then $\langle 1,16,16\rangle$ represents every positive $n$, except possibly when $a=0$ and $n\not\equiv1\pmod8$; when $a=1$ and $u\equiv2,3\pmod4$; or when $a\ge2$ and $u\equiv7\pmod8$. (See also the paragraph below.)
\item[{\bf (BF)}]\namedlabel{BF}{(BF)} (Brahmagupta--Fibonacci) $(a^2+b^2)(c^2+e^2)=(ac-be)^2+(ae+bc)^2$.
\item[{\bf (AM-GM)}]\namedlabel{AMGM}{(AM-GM)} (Inequality of arithmetic and geometric means) For real $a_1,\ldots,a_k\ge0$ we have $\frac{a_1+\cdots+a_k}{k}\ge(a_1\cdots a_k)^{1/k}$.
\end{itemize}

\textbf{On T11616:}  (From \cite{JP}.) The genus of $\langle 1,16,16\rangle$ (determinant $256$) has two classes; the other is
$4x^2+9y^2+9z^2+2yz+4xz+4xy$ (Jones--Pall's $(4,9,9,1,2,2)$), which is not regular: it misses exactly the squares $m^2$ with all prime factors of $m$ congruent to $1\pmod4$. Jones and Pall~\cite{JP} prove $\tern1{16}{16}$ to be regular (their Table~I lists $(1,r,r)$ for $r=16$ and Table~II the genus mate).

We also need one conjecture.

\begin{itemize}
\item[{\bf C145}]
(Conjecture 145.)\namedlabel{c:145}{C145}
$\langle1,4,5\rangle$ represents every positive $n$ except $n\equiv7\pmod8$, $n=4^a(8b+3)$ with $a\ge0$, and $n\in\{2,26,38\}$.
\end{itemize}

Conjecture C145 seems to be essential for the type of problem addressed in this article. We presume that it is well within reach of current theory, albeit nontrivial, possibly under the Generalized Riemann Hypothesis. (However, we ourselves are no experts on this.)
Results of Duke--Schulze-Pillot type~\cite{DSP} give it only for sufficiently large $n$ outside possible spinor-exceptional square classes, and ineffectively, so we treat it as a conjecture.
We verified C145 for $n\le3\cdot10^6$.

In fact all our conditional results only use $\tern145$ on numbers $m=4^au$ with $u\equiv6\pmod8$, so they follow from the weaker statement
\begin{itemize}
\item[{\bf C145$'$}]\namedlabel{c:145p}{C145$'$}
$x^2+4y^2+5z^2$ represents every positive integer $m=4^au$ with $u\equiv6\pmod 8$, except $m=38$.
\end{itemize} All our conditional results are therefore stated with C145$'$.

Statements \ref{T114} and \ref{T144} of the list at the start of this section are proved as follows.
\begin{lemma}\label{l:d}
$\tern114$ represents every positive $m\not\equiv3\pmod 4$ not of the form $4^a(8b+7)$; $\tern144$ represents every positive $m\equiv1\pmod4$ and every $m=4m'$ with $m'\ne4^a(8b+7)$.
\end{lemma}
\begin{proof}
These are both easy consequences of \ref{T111}. Since in both cases $m\ne4^a(8b+7)$ we may write $m=
p^2+q^2+r^2$. Now the number of odd elements among $p,q,r$ is congruent to $m$ modulo $4$, so:

Let $m\not\equiv3\pmod 4$. Then at least one of $p,q,r$ is even: say $r/2$ is an integer. Then 
 $m=p^2+q^2+4(r/2)^2$. 

Let $m\equiv1\pmod4$. Now two are even, so assume $q/2,r/2$ are integers. Then $m=p^2+4(q/2)^2+4(r/2)^2$.

Let  $m=4m'$ with $m'=p^2+q^2+r^2$. Then $m=(2p)^2+4q^2+4r^2$.
\end{proof}
\subsection{Some lemmas}\label{cov}

\begin{lemma}\label{l:types}
Every monomial square form $F=\sum_{i=1}^k M_i^{2}$ has a value-preserving specialization to a form whose variables have pairwise distinct fingerprints; in particular it has at most $2^k-1$ variables.
\end{lemma}

\begin{proof} If two variables $x,y$ of $F$ have the same fingerprint, then they only occur in $F$ through the product $xy$, so putting $y=1$ is a value-preserving specialization of $F$ with one variable fewer. Repeating this, we arrive at a form whose variables have pairwise distinct fingerprints. A fingerprint is a nonzero vector in $\{0,1\}^k$, so such a form has at most $2^k-1$ variables.
\end{proof} 

\begin{remark}\label{l:setone}
If $F\lp F'$, then every value of $F'$ is a value of $F$. In particular, if $F'$ is universal then so is $F$. 
\end{remark}

Since the forms will often specialize to several ternary forms, we want some tools that will help us to decide when a form is universal. 
Write $m=4^au$ with $4\nmid u$. For $a\ge3$, the integers $m$ not represented by the ternary forms listed in Section~\ref{Tools} fall into three cases: 
\begin{itemize}
\item[(A)] $u\equiv7\pmod 8$ for T111, T122, T114, T144, T11616; 
\item[(B)] $u\equiv14\pmod{16}$ for T112; 
\item[(C)] $u\equiv3\pmod 8$ for $\langle1,4,5\rangle$ under C145.
\end{itemize}

\begin{lemma}\label{l:cov}
Let $F$ have two branches $c+Q$ and $c+Q'$ with the same shift, where $Q,Q'$ are among the ternary forms above and belong to different cases among (A), (B), (C). Then every $n>c$ with $v_2(n-c)\ge6$ is represented by $F$ (using $\tern145$, and C145, only for $n-c>38$; if the cases are (B) and (C), C145$'$ suffices).
\end{lemma}
\begin{proof}
The large-valuation parts of the sets $E(Q)$ are: $v_2(m)$ even, $u\equiv7\pmod 8$ (A); $v_2(m)$ odd (B); $v_2(m)$ even, $u\equiv3\pmod 8$ (C): pairwise disjoint.
\end{proof}

\begin{example} \label{ex:1} $F=x^2+y^2+z^2+(xt)^2$. Specializing $t=0$ gives $\tern111$; $t=1$ gives $\tern211=\tern112$. Now $E(\tern111)\cap E(\tern112)=\emptyset$ so $F$ is universal.
\end{example}

\begin{lemma}[Scaling lemma]\label{l:scale}
Let $F=\sum_{i=1}^k M_i^2$ be a form in the variables $x_1,\dots,x_r$, and
suppose there are integers $e_1,\dots,e_r\ge 0$ such that, with $\deg(x_j)=e_j$,
all monomials have the same degree: $\deg(M_1)=\dots=\deg(M_k)=d$.
Then for every integer $c$
\[
F(c^{e_1}x_1,\dots,c^{e_r}x_r)=c^{2d}F(x_1,\dots,x_r).
\]
In particular, if $F$ represents $m$, then it represents $4^d m$. Hence $F$ is universal
as soon as it represents every $n\ge 1$ with $4^d\nmid n$.
\end{lemma}

\begin{example} \label{ex:2} (a class that needs the descent through scaling). 
\[ F=(1+y^2+z^2)x^2+t^2=x^2+(xy)^2+(xz)^2+t^2\]
 i.e.\ the class $\{x,xy,xz,t\}$ in the notation of Section~\ref{4less}, represents every $n\ge1$. Here $F(2x,y,z,2t)=4F$, so we may assume $n\not\equiv 0 \pmod 4$ (this is where we use the descent pointed out in Lemma~\ref{l:scale}). 
$F(1,y,z,t)=1+\langle 1,1,1\rangle$, so by \ref{T111} every $n$ not of the form $1+4^a(8b+7)$ is represented; it remains to treat $n=1+4^a(8b+7)$. If $a=0$ we get $8$ divides $n$, which is excluded, so we can assume $a\ge 1$. 
This means that $n\equiv 1 \pmod 4$. 
Now $n-4\equiv1\pmod4$, so $n-4$ is a value of $\tern144$ by \ref{T144}, and $n$ is a value of $F(2,y,z,t)=4+\tern144$. 
\end{example}

\begin{example} \label{ex:3} (conditional on C145$'$). $F=x^2+y^2+(xz)^2+(zt)^2$. $z=1$: $\tern{2}{1}{1}$; $z=2$: $\tern{5}{1}{4}=\tern145$. Cases (B) and (C) are disjoint, so every $n\ge1$ is represented (the values $2,26,38$ are not of the form $4^a(16b+14)$, so $\tern112$ covers them).
\end{example}

\subsection{The density exponent $e$}\label{expo}
This subsection describes a heuristic for how likely a form is to be universal. 
For a monomial square form $F=\sum_{i=1}^k M_i^2$ let
\[
N_F(N)=\#\bigl\{(|M_1(x)|,\ldots,|M_k(x)|):\ F(x)\le N\bigr\}
\]
be the number of tuples of monomial values with $F\le N$. We count tuples rather than points $x$: in $x^2+(xz)^2+\dots$ with $x=0$ every $z$ gives the same tuple, whereas there would be infinitely many points. Writing $r_F(n)$ for the number of tuples with $M_1^2+\dots+M_k^2=n$, we have $N_F(N)=\sum_{n\le N}r_F(n)$, so $N_F(N)/N$ is the average number of representations of an integer $n\le N$.

To estimate $N_F(N)$, first count the tuples with all $M_i\ne0$. Give each variable a weight $w_j\ge0$ and think of $|x_j|\approx N^{w_j}$. If every monomial has weighted degree at most $\tfrac12$, then $F(x)$ is at most about $N$, and there are about $N^{\sum_j w_j}$ such points. Taking the best weights suggests that there are about $N^{e_F}$ such tuples, where
\[
e_F=\max\Bigl\{\textstyle\sum_j w_j:\ w\ge0,\ \sum_{j\in M}w_j\le\tfrac12\ \text{for every monomial }M\Bigr\}.
\]
However, this only counts the tuples with all $M_i\ne0$. A tuple in which some monomials are $0$ comes from a specialization $F'$ of $F$ that puts some variables equal to $0$ (written in the variables that still occur in it), and there may be many more of those. So we expect $N_F(N)=N^{e+o(1)}$ with
\[ e=\max_{F'} e_{F'}, \]
where $F'$ runs over $F$ itself and these specializations. For example, for $F=x^2+y^2+z^2+t^2+(xyztu)^2$ we get $e_F=\tfrac12$, because nonzero $x,y,z,t,u$ force $|xyztu|\le\sqrt N$; but $u=0$ leaves $x^2+y^2+z^2+t^2$, with about $N^2$ tuples, so $e=2$.

The usual refinement, from counting lattice points in a dilated polytope, is $N_F(N)\asymp N^e(\log N)^d$, where $d$ is the dimension of the face of optimal weightings, taken for the form $F'$ attaining $e$ (the largest one if several forms attain $e$).

For the five-monomial classes of Section~\ref{5m5v} that do not miss a small value, $e$ is $1$, $\tfrac32$, $2$ or $\tfrac52$, and the open ones have $e=1$ or $e=\tfrac32$. For $e=\tfrac32$ the average number of representations grows like $\sqrt N$. For $e=1$ it grows like $(\log N)^d$; it would stay bounded only if $d=0$, and only then would one expect infinitely many missed values. All twelve five-variable classes with $e=1$ have $d=2$ or $d=3$ (in each of them the form itself attains the largest $d$), yet two of them are not universal (Section~\ref{s:5over}).

These heuristics prove {\em nothing:} many non-universal forms have $e=1$ or even $e=\tfrac32$; for example $\langle 1,1,1\rangle$ has $e=\tfrac32$ and misses all $4^a(8b+7)$. They only indicate where to place one's bets.

\section{Forms with at most four monomials}

\label{4less}

{\bf At most three monomials:}
 A sum of at most three squares misses $7$, so no form with $k\le3$ is universal, in any number of variables.
So this section is actually on the case of exactly four monomials.

\subsection{The 1885 renaming classes}\label{s:4types}

Let $F=M_1^2+M_2^2+M_3^2+M_4^2$ be a monomial square form with four distinct nonconstant monomials. By Lemma~\ref{l:types}, $F$
can be assumed to have at most $2^4-1=15$ variables.
Counting these forms up to renaming (Appendix~\ref{app:4}) gives {\bf 1885} renaming classes of such four-monomial forms.
These are the objects of this section, and they have between $3$ and $15$ variables.

We denote a class by the set of its monomials, e.g.\ $\{xt,\,y,\,z,\,t\}$ means $x^2t^2+y^2+z^2+t^2$. Classes with at most six variables use the letters $x,y,z,t,u,v$.

By an exact computation of all values up to $10^4$ (Appendix~\ref{app:4}) the $1885$ classes fall into three groups.
\begin{itemize}
\item $944$ specialize to one of the forms \#1--\#4 of Table~\ref{tab:four} and are universal ($300$ of them only under C145$'$, see Theorem~\ref{t:four}).
\item $841$ miss a value $\le143$.
\item  $100$ classes, with $5$ to $12$ variables, all have $9591$ as smallest missed value. They are exactly the classes that specialize to
\[
G=\{xzu,\ yt,\ zu,\ tu\}
\]
and none of \#1--\#4. The form $G$ has five variables.
\end{itemize}

\subsection{The minimal universal forms}
\begin{center}\small
\captionof{table}{The four minimal universal four-monomial forms}\label{tab:four}
\begin{tabular}{llcp{7.4cm}}
\toprule
\# & monomials & variables & status\\
\midrule
1 & $\{x,\,y,\,z,\,t\}$ & 4 & universal (\ref{T1111})\\
2 & $\{xt,\,y,\,z,\,t\}$ & 4 & universal (Example~\ref{ex:1}, using \ref{T111} and \ref{T112})\\
3 & $\{xt,\,yt,\,z,\,t\}$ & 4 & universal (Example~\ref{ex:2}, using \ref{T111}, \ref{T144} and Lemma~\ref{l:scale})\\
4 & $\{xz,\,y,\,zt,\,t\}$ & 4 & \leavevmode universal if C145$'$ holds (Example~\ref{ex:3})\\
\bottomrule
\end{tabular}
\end{center}

\begin{theorem}\label{t:four}
The forms \#1--\#4 are minimal, and 
a four-monomial form is (conditionally) universal if and only if it has a renamed specialization to 
 one of \#1--\#4, where \#4 gives universality only under C145$'$. 
\end{theorem}
\begin{proof}
A brief analysis shows that every specialization of \#1--\#4 that is not a renaming misses a small value that the form itself represents (for \#4 this is a direct check, not using C145$'$), so these four forms are minimal. 
Appendix~\ref{app:4} shows that either a form is universal and specializes to one of \#1--\#4, or it misses as smallest value  one of $2,3,6,13,14,15,29,79,143$ or $9591$.
\end{proof}

In fact, the list of smallest missed values in this proof yields a version of the ``15-theorem'' for these forms:

\begin{theorem}\label{t:fifteen}
A four-monomial square form is universal if and only if it represents  14, 15, 29 and 9591 (for a form that specializes to \#4 but to none of \#1--\#3, the ``if'' part assumes C145$'$). 
\end{theorem}

The proof is perhaps a bit unsatisfying, as it simply follows from the classification, together with the fact that each of the four numbers is the smallest missed value of some form. The set is not unique: $14$, $15$ and $29$ can be replaced by other values, but they are the smallest possible choices.

\subsection{The 100 classes missing 9591 and their king $H$}\label{s:9591}

An interesting curiosity is that among these $100$ classes  there is a largest one, and it specializes to all the others. 
It is the form  
\[H=\{x_1x_3x_4x_6x_7x_9x_{10}x_{12}, x_2x_3x_5x_6x_8x_9x_{11}x_{12},
x_4x_5x_6x_{10}x_{11}x_{12}, x_7x_8x_9x_{10}x_{11}x_{12} \}. 
\]
Every variable occurs in the first or the second monomial.
Setting all variables except $x_1,x_2,x_4,x_8,x_{10}$ equal to $1$ gives $
\{x_1x_4x_{10}, x_2x_8, x_4x_{10},x_8x_{10}\}$, which is a renaming of $G$.

\begin{proposition}\label{p:ruled}
The $100$ classes are exactly the classes $F$ with $H\lp F\lp G$. In particular $E(G)\supseteq E(F)\supseteq E(H)\supseteq \{9591\}$.
\end{proposition}
\begin{proof}
It follows essentially from the classification. \end{proof}

\emph{Is $9591$ the only value missed by $H$?} Only integers $n=4^a(8b+7)$ need checking, since $H$ specializes to $\tern111$ ($x_1=0$, $x_3=x_5=x_6=x_8=x_9=x_{10}=x_{11}=x_{12}=1$) and \ref{T111} applies. A check up to $10^7$ shows no other missed values. 

\section{Five monomials in at most five variables}

\label{5m5v}

\subsection{Setting}
We will use the variables $x,y,z,t,u$ in this setting. 
As in Section~\ref{4less}, we write a form as the set of its monomials: for example
\[ \{x,xy,yz,t,tu\} \textup{~means~} x^2+(xy)^2+(yz)^2+t^2+(tu)^2.\]
There are $31$ squarefree monomials $\ne1$ in five variables, giving $\binom{31}{5}=169\,911$ forms, which fall into $2160$ renaming classes. 

\subsection{Why stop at five variables}

When looking at Lemma~\ref{l:types} we see that we only have a bound of $2^5-1=31$ variables, and 5 is far off. However, already 6 variables instead of 5 means that the search space gets enormous (41 times as big as the 5-variable case). Also, the number of open cases with five variables is already large enough that one may wonder what six variables would yield beyond a lot of data.

\subsection{Overview of the results}\label{s:5over}
We have 2160 classes. Of these, 2077 are not universal. The 83 remaining ones fall into the following categories:
\begin{center}
\captionof{table}{Status of the $83$ remaining classes}\label{tab:83}
\begin{tabular}{lrl}
\toprule
status & classes & forms\\
\midrule
proved by the 15-theorem (Table~\ref{tab:settled}, rows 1--11) & 11 & 271\\
branches, unconditional (Table~\ref{tab:settled}, rows 12--48) & 37 & 3\,150\\
branches, if C145$'$ holds (Table~\ref{tab:settled}, rows 49--66) & 18 & 2\,100\\
open (Table~\ref{tab:open}) & 17 & 1\,860\\
\midrule
total & 83 & 7\,381\\
\bottomrule
\end{tabular}
\end{center}

The $17$ open classes are listed in Table~\ref{tab:open}.

(Rows 1--7 are checked up to $10^7$ and rows 8--17 up to $10^9$; rows 8--14, 16 and 17 are also checked up to $10^{10}$ in the $60$ residue classes modulo $720$ where they have the fewest representations.)
\small
\begin{longtable}{r p{6cm} c l}
\caption{The $17$ open classes}\label{tab:open}\\
\toprule
\# & $\{M_i\}$ & exponent $e$ & ternary branch available\\
\midrule\endfirsthead
\toprule
\# & $\{M_i\}$ & exponent $e$ & ternary branch available\\
\midrule\endhead
1 & $\{x,\,y,\,xz,\,yt,\,ztu\}$ & $\tfrac32$ & yes (6 branches) \\
2 & $\{x,\,y,\,xz,\,xzt,\,ztu\}$ & $\tfrac32$ & yes (43 branches) \\
3 & $\{x,\,xy,\,xz,\,xt,\,yzu\}$ & $\tfrac32$ & yes (31 branches) \\
4 & $\{x,\,xy,\,xz,\,yt,\,zu\}$ & $\tfrac32$ & yes (39 branches) \\
5 & $\{x,\,xy,\,xz,\,yt,\,yzu\}$ & $\tfrac32$ & yes (40 branches) \\
6 & $\{x,\,xy,\,yz,\,tu,\,xyt\}$ & $\tfrac32$ & yes (43 branches) \\
7 & $\{x,\,xy,\,yz,\,xyt,\,ytu\}$ & $\tfrac32$ & yes (39 branches) \\
8 & $\{x,\,y,\,xz,\,yt,\,xzu\}$ & $1$ & no \\
9 & $\{x,\,y,\,xz,\,xyt,\,xzu\}$ & $1$ & no \\
10 & $\{x,\,y,\,xz,\,xzt,\,yzu\}$ & $1$ & no \\
11 & $\{x,\,y,\,xz,\,xzt,\,ytu\}$ & $1$ & no \\
12 & $\{x,\,y,\,xz,\,xzt,\,xztu\}$ & $1$ & no \\
13 & $\{x,\,xy,\,xz,\,yt,\,xzu\}$ & $1$ & no \\
14 & $\{x,\,xy,\,xz,\,yt,\,ytu\}$ & $1$ & no \\
15 & $\{x,\,xy,\,xz,\,yt,\,ztu\}$ & $1$ & no \\
16 & $\{x,\,xy,\,yz,\,zu,\,xyt\}$ & $1$ & no \\
17 & $\{x,\,xy,\,yz,\,xyt,\,yzu\}$ & $1$ & no \\
\bottomrule
\end{longtable}
\normalsize

The $11$ classes settled by the 15-theorem have density exponent $e=2$, except $\{x,y,z,t,u\}$ with $e=\tfrac52$; the other settled classes have $e=\tfrac32$. (Their certificates are in Table~\ref{tab:settled} in Appendix~\ref{app:5}.)
The ten open classes with $e=1$ (rows 8--17) have far fewer representations per integer than the others, but their averages still grow like $(\log n)^2$ or $(\log n)^3$. Two further classes with $e=1$ turn out to miss rather large values: 

\begin{center}\small
\captionof{table}{The two $e=1$ classes with a large missed value}\label{tab:e1miss}
\begin{tabular}{ll}
\toprule
class & missed values\\
\midrule
$\{x,\,xy,\,xz,\,yt,\,xytu\}$ & $274647,\ 3885279$\\
$\{x,\,xy,\,xz,\,yt,\,xyzu\}$ & $648798,\ 2464638,\ 3550398,\ 26671358,\ 178039518,\ 390944958$\\
\bottomrule
\end{tabular}
\end{center}
For both classes these are all missed values up to $10^9$. For the second class, the search up to $10^{10}$ in the residue classes mentioned with Table~\ref{tab:open} also found $1445025918$.

It is quite possible that more can be eliminated, but our search did not find any further missed values (see the search bounds given with Table~\ref{tab:open}).

A ternary branch alone gives of the order $N^{3/2}$ points with $F\leq N$, which partly explains why the $e=1$ cases are exactly those without ternary branches.

\subsubsection*{Why the open cases remain open}
Rows 2, 4, 5, 6 and 7 of Table~\ref{tab:open} 
fail to be settled affirmatively for the \emph{same} reason. Each has the descent of Lemma~\ref{l:scale} (weights $(1,0,0,1,1)$, $(1,1,0,0,1)$ or $(1,0,1,0,1)$, all monomials of degree $1$), so a smallest counterexample $n$ is squarefree. Fixing two variables $y=b$, $z=c$ (resp.\ $y=b,t=d$; for row 2, $z=b,t=d$) leaves a diagonal ternary form $\langle 1+b^2+c^2,\,b^2,\,c^2\rangle$, $\langle1+b^2+c^2,b^2,b^2c^2\rangle$, $\langle 1+b^2+b^2d^2,b^2,d^2\rangle$, $\langle 1+b^2+b^2d^2,1,b^2d^2\rangle$ or $\langle 1+b^2+b^2d^2,b^2,b^2d^2\rangle$ (for rows 4, 5, 6, 2 and 7, respectively); for $b=1$ and $c$ (resp.\ $d$) equal to $1$ or $2$ these families contain the same two forms
\[
\begin{array}{ll}
\tern113:& \text{misses } 9^k(9m+6),\qquad\\
 \tern146:& \text{misses } n\equiv2\pmod{16}\ \text{and}\ n=3^{\mathrm{odd}}u,\ u\equiv1\pmod 3.
\end{array}
\]
Their exceptional sets overlap exactly in the squarefree numbers
\[
n\equiv 2\pmod{16},\ n\equiv 6\pmod 9\ \Longleftrightarrow\ n\equiv114\pmod{144}.
\]
 Each of these five rows would follow from showing that the form represents every squarefree $n\equiv114\pmod{144}$. 

Row~1 has different obstructions: the branches with shift $0$ lead to $\tern122$ and leave $n\equiv7\pmod 8$, which would need the non-regular forms $\tern245$ and $\tern{5}{10}{36}$ (for the regular ternary forms see~\cite{JKS}).

Of the $31$ ternary branches of row~3, only $1+\tern111$ turned out to be useful.

\subsection{Finitely many test values}\label{s:test5}
As for four monomials (Theorem~\ref{t:fifteen}), the computations of this section give a finite set of test values, but here only conditionally.
\begin{proposition}\label{p:test5}
Let \[S=\{15,\,30,\,78,\,2159,\,3423,\,7280,\,20959,\,39759,\,63985,\,65663,\,274647,\,648798\}.\] A five-monomial form in five variables that represents every element of $S$ belongs to one of the $83$ classes of Section~\ref{s:5over}. Hence, if C145$'$ holds and the $17$ classes of Table~\ref{tab:open} are universal, a five-monomial form in five variables is universal if and only if it represents every element of $S$.
\end{proposition}
\begin{proof}
Each of the $2077$ non-universal classes misses an element of $S$; this is an exact finite computation (of all values up to $10^5$ of every class, see Appendix~\ref{app:5}, and a direct check of $274647$ and $648798$ for all classes; one class, $\{x,y,xz,xzt,xyzu\}$, misses no element of $S$ below $10^5$ and is caught only by $648798$). The second statement then follows from Table~\ref{tab:settled} in Appendix~\ref{app:5}.
\end{proof}
The twelve classes below miss pairwise disjoint sets of integers in $[1,10^5]\cup\{274647,648798\}$, so no set of fewer than twelve such integers works. (We did not compute all misses above $10^5$, so a larger integer missed by two of them at once is not excluded.)
\begin{center}\small
\captionof{table}{Twelve classes with pairwise disjoint missed values}\label{tab:twelve}
\begin{tabular}{rlrlr}
\toprule
\# & class & misses $\le10^5$ & smallest misses & in $S$\\
\midrule
1 & $\{x,\,xy,\,yz,\,zt,\,tu\}$ & 3 & 15, 31, 255 & 15\\
2 & $\{x,\,xy,\,yz,\,yt,\,xzu\}$ & 2 & 30, 3198 & 30\\
3 & $\{x,\,xy,\,yz,\,xyt,\,ztu\}$ & 134 & 78, 142, 191, \dots & 78\\
4 & $\{x,\,y,\,xz,\,xyt,\,ztu\}$ & 19 & 79, 240, 1919, \dots & 3423\\
5 & $\{x,\,xy,\,xz,\,xt,\,xyzu\}$ & 851 & 113, 241, 449, \dots & 63985\\
6 & $\{x,\,y,\,xz,\,yt,\,xyzu\}$ & 2 & 519, 20959 & 20959\\
7 & $\{x,\,y,\,xz,\,xzt,\,xytu\}$ & 4 & 2159, 30239, 30623, 68543 & 2159\\
8 & $\{x,\,y,\,z,\,xyt,\,xzu\}$ & 3 & 7280, 18880, 68080 & 7280\\
9 & $\{x,\,xy,\,xz,\,xyt,\,ytu\}$ & 1 & 39759 & 39759\\
10 & $\{x,\,xy,\,xz,\,yt,\,yztu\}$ & 3 & 65663, 66237, 89183 & 65663\\
11 & $\{x,\,xy,\,xz,\,yt,\,xytu\}$ & 0 & 274647, \dots & 274647\\
12 & $\{x,\,xy,\,xz,\,yt,\,xyzu\}$ & 0 & 648798, \dots & 648798\\
\bottomrule
\end{tabular}
\end{center}
The set $S$ is not unique. Among the sets of twelve integers from $[1,10^5]\cup\{274647,648798\}$ that work, $S$ is the one whose largest element is as small as possible, then its second largest element, and so on. The elements $39759$, $274647$ and $648798$ are forced: rows 9, 11 and 12 miss no other integer up to $10^5$, $10^6$ and $10^6$, respectively.

\section{The $P_{m,d}$ and $C_{m,d}$ forms}

\label{PnmCnm}

\subsection{Setting and main results}
For integers $m\ge2$ and $d\ge2$ put $n=m+d-1$ and let
\begin{gather*}P_{m,d}(x)=\sum_{i=1}^{d}(x_ix_{i+1}\cdots x_{i+m-1})^2,\\
C_{m,d}(x)=\sum_{i=1}^{n}(x_ix_{i+1}\cdots x_{i+m-1})^2\quad(x_{i+n}=x_i).\end{gather*}
So $P_{m,d}$ has $d$ monomials and $C_{m,d}$ has $n$ monomials, each a product of $m$ consecutive variables, and both have $n$ variables; $C_{m,d}$ is $P_{m,d}$ plus $m-1$ more monomials.
For $d=1$ we only use the single square $P_{m,1}=(x_1\cdots x_m)^2$, in Lemma~\ref{lem:zero}(a),(b) for $d=2$.

  For $C_{m,d}$ it is useful to split off the values at points with all $x_i\ne0$:
\[ \mathcal{C}^*_{m,d}=\{C_{m,d}(x):\ x\in\N^n,\ x_1\cdots x_n\ne0\}\subseteq[n,\infty), \]
since each of the $n$ monomials is then at least $1$.

Tables~\ref{tab:Pmd} and~\ref{tab:Cmd} list the results of this section: they give the least positive integer that is not represented by the form. 
``U'' means universal (proved), ``U$'$'' means universal if C145$'$ holds, and ``$\infty$'' means that no positive integer up to $10^7$ is missed, but universality is open.

\begin{center}\small\renewcommand{\arraystretch}{1.2}
\begin{minipage}{\linewidth}\centering
\captionof{table}{Least positive integer not represented by $P_{m,d}$}\label{tab:Pmd}
\begin{tabular}{@{}lccccc@{}}
\toprule
$m\downarrow,\ d\rightarrow$ & $2$ & $3$ & $4$ & $5\le d\le m+3$ & $d\ge m+4$\\\midrule
$2$ & $3$ & $7$ & U$'$ & U & U\\
$3$ & $3$ & $7$ & $9591$ & $\infty$ & U\\
$\ge4$ & $3$ & $7$ & $9591$ & $\infty$ & U\\\bottomrule
\end{tabular}
\end{minipage}

\bigskip
\begin{minipage}{\linewidth}\centering
\captionof{table}{Least positive integer not represented by $C_{m,d}$}\label{tab:Cmd}
\begin{tabular}{@{}lcccccc@{}}
\toprule
$m\downarrow,\ d\rightarrow$ & $2$ & $3$ & $4$ & $5$ & $6\le d\le m+4$ & $d\ge m+5$\\\midrule
$2$ & $2$ & $3$ & $7$ & U$'$ & U & U\\
$3$ & $2$ & $3$ & $7$ & $\infty$ & $\infty$ & U\\
$4$ & $2$ & $3$ & $15$ & $9591$ & $\infty$ & U\\
$\ge5$ & $2$ & $3$ & $7$ & $9591^{*}$ & $\infty$ & U\\\bottomrule
\end{tabular}
\end{minipage}
\end{center}
${}^*$ For every $m\ge5$ with $m\ne9587$ the smallest missed value is $9591$; for $m=9587$ it is $77871$ (Proposition~\ref{p:C5}).

\noindent The entries $\infty$ for $P_{m,d}$ ($m\ge3$, $5\le d\le m+3$) and for $C_{m,d}$ ($m\ge3$, $6\le d\le m+4$) all become U if row 8 or row 14 of Table~\ref{tab:open} (Section~\ref{5m5v}) is universal (Proposition~\ref{p:PtableA}). The entry $\infty$ for $C_{3,5}$ does not reduce to Section~\ref{5m5v} in this way.

The remainder of this section is dedicated to proving the entries of Tables~\ref{tab:Pmd} and~\ref{tab:Cmd}.

\subsection{Basic lemmas}
\begin{lemma}\label{lem:zero}
Let $n=m+d-1$.
\begin{enumerate}
\item[(a)] If $d\ge2$, then $P_{m,d-1}(\N^{n-1})\subseteq  P_{m,d}(\N^n)$.
\item[(b)] If $d\ge2$, then $C_{m,d}(\N^n)=P_{m,d-1}(\N^{n-1})\cup \mathcal{C}^*_{m,d}$.
\item[(c)] $\mathcal{C}^*_{m,d}\subseteq[n,\infty)$ and $n\in \mathcal{C}^*_{m,d}$.
\end{enumerate}
\end{lemma}
\begin{proof}
(a) $(P_{m,d})_{x_n=0}=P_{m,d-1}$. 
(b) We have $(C_{m,d})_{x_n=0}=P_{m,d-1}$, and by the cyclic symmetry the same holds for $x_i=0$ up to renaming. So the values at points with some $x_i=0$ are exactly $P_{m,d-1}(\N^{n-1})$; the other values form $\mathcal{C}^*_{m,d}$.
(c) The lowest value of $\mathcal{C}^*_{m,d}$ is attained when all $x_i=1$, and $C_{m,d}(1,\ldots,1)=n$.
\end{proof}

\begin{remark}\label{lem:core}
Let $m\ge d\ge3$. 
Then 
\[P_{m,d}(x)=(x_d\cdots x_m)^2\,P_{d-1,d}(x_1,\ldots,x_{d-1},x_{m+1},\ldots,x_{m+d-1}).\]
Hence 
\[P_{m,d}(\N^n)=\bigcup_{\lambda\in\N} \lambda^2P_{d-1,d}(\N^{2d-2}).\]
In particular, the set of values does not depend on $m$ (for $m\ge d$), and a squarefree integer is a value of $P_{m,d}$ if and only if it is a value of $P_{d-1,d}$.
\end{remark}

\subsection{The table for $P_{m,d}$}\label{s:Ptab}

\paragraph{\boldmath $d=2,3$:} Here $P_{m,d}$ is a sum of two or three squares, and the table entries follow by a short computation.

\paragraph{\boldmath $d=4$, $m=2$:} Setting $x_2=1$ in $P_{2,4}=(x_1x_2)^2+(x_2x_3)^2+(x_3x_4)^2+(x_4x_5)^2$ gives $x_1^2+x_3^2+(x_3x_4)^2+(x_4x_5)^2$, which after renaming $(x_1,x_3,x_4,x_5)\to(y,t,z,x)$ is class \#4  of Table~\ref{tab:four} in Section~\ref{4less}, $\{xz,\,y,\,zt,\,t\}$. So $P_{2,4}$ is universal if C145$'$ holds (Theorem~\ref{t:four}).

\paragraph{\boldmath $d=4$, $m\geq 3$:}
For $m=3$, setting $x_4=1$ in $P_{3,4}$ gives $\{x_1x_2x_3,\,x_2x_3,\,x_3x_5,\,x_5x_6\}$, which is $G=\{xzu,yt,zu,tu\}$ from Section~\ref{s:9591} after renaming $(x_1,x_2,x_3,x_5,x_6)\to(x,z,u,t,y)$. The form $P_{3,4}$  is one of the $100$ classes of Section~\ref{s:9591}, and its smallest missed value is $9591$. For $m\ge4$, Remark~\ref{lem:core} gives 
$P_{m,4}(\N^n)=\cup_{\lambda \in \N} \lambda^2 P_{3,4}(\N^6)$. Since $9591=3\cdot23\cdot139$ is squarefree, the smallest missed value of $P_{m,4}$ is also $9591$.
(In fact, up to $10^7$ the only values missed by $P_{m,4}$, $m\geq 3$, are $9591$ and $77871=3\cdot 101\cdot 257$: for $m=3$ by a finite computation, and for $m\ge4$ by the same argument as above.)

\paragraph{\boldmath $d\ge5$:} This takes more work.
\begin{theorem}\label{thm:univ}
\begin{enumerate}
\item[(a)] $P_{2,d}$ is universal for every $d\ge5$.
\item[(b)] For $m\ge2$ and $d\ge m+4$, $P_{m,d}$ is universal.
\end{enumerate}
\end{theorem}
\begin{proof}
(a) We specialize $P_{2,5}=(x_1x_2)^2+(x_2x_3)^2+(x_3x_4)^2+(x_4x_5)^2+(x_5x_6)^2$ in two ways: $x_2=1,x_3=0,x_5=1$ gives $\tern111$, $x_1=0,x_3=1,x_5=1$ gives $\tern112$. 
$E(\tern111)\cap E(\tern112)=\emptyset$ so $P_{2,5}$ is universal.
Universality of $P_{2,d}$ now follows from Lemma~\ref{lem:zero}(a).  (With $x_5=1$, $P_{2,5}$ is row 34 of Table~\ref{tab:settled}.)

(b) By Lemma~\ref{lem:zero}(a) it suffices to treat $d=m+4$, so $n=2m+3$. The variable $x_{m+2}$ lies in the monomials $i$ with $3\le i\le m+2$. Putting $x_{m+2}=0$ leaves only four monomials nonzero. Additionally setting all remaining variables except 
$x_1,x_{m+1},x_{m+3},x_{2m+3}$ to 1 we obtain the form 
$x_1^2+x_{m+1}^2+x_{m+3}^2+x_{2m+3}^2$, which is universal by \ref{T1111}.
\end{proof}

It remains to justify the entries $\infty$, for $m\ge3$ and $5\le d\le m+3$. By Lemma~\ref{lem:zero}(a), $P_{m,d}$ has all values of $P_{m,4}$, so by the case $d=4$ it suffices that $9591$ and $77871$ are values of $P_{m,5}$. Indeed
\begin{gather*}P_{3,5}(1,3,17,1,1,2,32)=9591,\quad P_{3,5}(1,3,53,1,1,70,2)=77871,\\
P_{4,5}(1,3,1,7,1,1,2,46)=9591,\quad P_{4,5}(1,3,1,13,1,1,2,136)=77871,\end{gather*}
and for $m\ge5$, $P_{m,5}$ has all values of $P_{4,5}$ by Remark~\ref{lem:core}.

Whether these forms are universal comes down to Section~\ref{5m5v}:
\begin{proposition}\label{p:PtableA}
Let $m\ge3$. Putting all variables of $P_{m,5}$ except $x_1,x_2$ and $x_{m+2},x_{m+3},x_{m+4}$ equal to $1$ gives row 8 of Table~\ref{tab:open}, and putting all except $x_1,x_3$ and $x_{m+2},x_{m+3},x_{m+4}$ equal to $1$ gives row 14. So if row 8 or row 14 of Table~\ref{tab:open} is universal, then $P_{m,d}$ is universal for all $m\ge3$ and $d\ge5$, and $C_{m,d}$ for all $m\ge3$ and $d\ge6$.
\end{proposition}
\begin{proof}
The five monomials become $x_1x_2$, $x_2$, $x_{m+2}$, $x_{m+2}x_{m+3}$, $x_{m+2}x_{m+3}x_{m+4}$, which is $\{x,y,xz,yt,xzu\}$ (row 8) for $(x,y,z,t,u)=(x_{m+2},x_2,x_{m+3},x_1,x_{m+4})$. In the second case they become $x_1x_3$, $x_3$, $x_3x_{m+2}$, $x_{m+2}x_{m+3}$, $x_{m+2}x_{m+3}x_{m+4}$, which is $\{x,xy,xz,yt,ytu\}$ (row 14) for $(x,y,z,t,u)=(x_3,x_{m+2},x_1,x_{m+3},x_{m+4})$. The rest is Lemma~\ref{lem:zero}(a),(b).
\end{proof}
\subsection{The table for $C_{m,d}$}\label{s:Ctab}
By Lemma~\ref{lem:zero}(b),(c), the values of $C_{m,d}$ below $n$ are exactly those of $P_{m,d-1}$, and $n=C_{m,d}(1,\ldots,1)$ is a value.

\paragraph{\boldmath $2\le d\le4$:} These cases are straightforward computations. We treat $m=d=4$ as an example, since it is the exception in its column ($15$ instead of $7$). 
By Lemma~\ref{lem:zero}(b)
we get $C_{4,4}(\N^7)=P_{4,3}(\N^6)\cup \mathcal{C}^*_{4,4}$. 
 The values $P_{4,3}(u,1,1,1,w,0)=u^2+w^2$, $P_{4,3}(u,1,1,1,1,w)=u^2+1+w^2$ and $P_{4,3}(1,1,2,1,1,1)=12$ give $1,\ldots,6,8,\ldots,14$. 

Since $15\equiv7\pmod 8$, it is not a value of $P_{4,3}$ (which is a sum of three squares).
Now we see that $\mathcal{C}^*_{4,4}$ has $7$ and $19$ as its two smallest elements
since $C_{4,4}(1,1,1,1,1,1,1)=7$ and $C_{4,4}(1,1,1,1,1,1,2)=19$ (if some $x_i\ge2$, the $4$ monomials containing $x_i$ give at least $4\cdot 4$, the other $3$ at least $3$).
Hence $7$ is a value, and $15$ is the smallest missed value of $C_{4,4}$.

\paragraph{\boldmath $d\ge5$, except ($d=5$, $m\ge4$):} By Lemma~\ref{lem:zero}(b), $C_{m,d}(\N^n)\supseteq P_{m,d-1}(\N^{n-1})$, so Table~\ref{tab:Pmd} gives the entries U$'$ for $C_{2,5}$, U for $C_{2,d}$ with $d\ge6$ and for $d\ge m+5$, and $\infty$ for $m\ge3$, $6\le d\le m+4$. For $C_{3,5}$, the only values up to $10^7$ that $P_{3,4}$ misses are $9591$ and $77871$ (Section~\ref{s:Ptab}), and
\[C_{3,5}(1,1,4,14,1,3,9)=9591,\qquad C_{3,5}(1,1,3,9,1,26,4)=77871,\]
so $C_{3,5}$ represents every positive integer up to $10^7$.

\paragraph{\boldmath $d=5$, $m\ge4$:}
\begin{proposition}\label{p:C5}
Let $m\ge4$. If $m\ne9587$, the least positive integer not represented by $C_{m,5}$ is $9591$. If $m=9587$, it is $77871$.
\end{proposition}

For most other results we have omitted the details of the computations, which are feasible after some obvious reductions. For Proposition~\ref{p:C5} we include them, since they are essential to the proof, not obvious, and in our opinion interesting in their own right.

\begin{proof}(of Proposition~\ref{p:C5}.)
Recall $n=m+4$, the number of variables.
A priori the problem is infinite, since the number of monomials grows with $m$ (unlike for $P_{m,d}$, cf.\ Remark~\ref{lem:core}). We first bound $m$.

Lemma~\ref{lem:zero}(b) gives $C_{m,5}(\N^n)=P_{m,4}(\N^{n-1})\cup \mathcal{C}_{m,5}^*$.
 $P_{m,4}$ represents every positive integer below $9591$, and its only missed values up to $10^7$ are $9591$ and $77871$ (Section~\ref{s:Ptab}).
It remains to decide whether these two numbers lie in $\mathcal{C}_{m,5}^*$. We will consider the lowest values appearing in $\mathcal{C}^*_{m,5}$.
Note that all values in $\mathcal{C}^*_{m,5}$ are obtained by substituting $x_i\geq 1$ only.

{\em Substituting at most one value $\ne1$:} Say $x_j=a\geq 1$, then this yields $4+ma^2$ (giving $4+m, 4+4m, 4+9m,\ldots$).

{\em Substituting at least two values $x_j, x_k\geq 2$:} 
There are $n=m+4$ monomials, each a product of $m$ cyclically consecutive variables. So each monomial leaves out exactly 4 variables, and each variable is left out of exactly 4 monomials. This means that at most 8 monomials leave out one of $x_j$ or $x_k$. Thus, at least $n-8=m-4$ monomials contain both (note $m-4\geq 0$). 
So the value of $C_{m,5}$ is at least $(n-8)(x_jx_k)^2+8\cdot 1^2\geq (m-4)(2\cdot 2)^2+8=16(m-4)+8$.

Now: if $m\geq 603$, then $16(m-4)+8\geq 16(603-4)+8=9592>9591$ so it can only be obtained by $4+ma^2=9591$, or $ma^2=9587$. This number is prime, however, so $a=1, m=9587$ is the only solution. So for $m\ge603$ we have $9591\in\mathcal{C}^*_{m,5}$ only for $m=9587$. 
For $m=9587$ the values of $\mathcal{C}^*_{9587,5}$ are $9591$, the numbers $4+9587a^2$ with $a\ge2$, which are $38\,352,\,86\,287,\ldots$, and numbers at least $16\cdot9583+8=153\,336$. So $77871\notin\mathcal{C}^*_{9587,5}$.

For $4\le m\le602$, deciding whether $9591\in\mathcal{C}^*_{m,5}$ is a finite search, but a very large one. The following argument makes it feasible.
 Define $X=x_1\cdots x_n$. Since each $x_i$ appears in $m$ of the monomials of 
 $C_{m,5}=\sum_{i=1}^n M_i^2$ where $M_i=x_ix_{i+1}\cdots x_{i+m-1}$ we see that the product $M_1M_2\cdots M_n$ of all these monomials equals $X^m$.  This allows us to use \ref{AMGM}, as $C_{m,5}$ is a sum which can be bounded from below:
\[ C_{m,5}=\sum_{i=1}^n M_i^2\geq n (M_1^2M_2^2\cdots M_n^2)^{\frac{1}{n}}
= nX^{\frac{2m}{n}}.\]
Thus, if $C_{m,5}(x)=9591$, then $9591\ge n\,X^{2m/n}$, so $X\le(9591/n)^{n/(2m)}$; for $m\ge12$ we have $n\ge16$ and $2m/n\ge\frac32$, so $9591\ge16X^{3/2}$, which gives $X\le71$.

 For $4\le m\le 11$, the bound $X\le(9591/n)^{n/(2m)}$ gives $X\le 1198$, $530$, $305$, $204$, $150$, $117$, $96$, $81$ for $m=4,\ldots,11$. A point $x\in(\N^*)^n$ is determined by the positions and the values of its coordinates that are at least $2$. For these $m$ there are between $357\,704$ (for $m=11$) and $3\,384\,286$ (for $m=4$) such points, about $8.3\cdot10^6$ in total, and none of them gives $9591$.
 
For $12\le m\le 602$ we have $X\le 71$. Listing all points is no longer feasible (for $m=602$ there are about $7.5\cdot10^{13}$ of them), but it is not necessary. 
First, notice that due to $X\le 71$ at most 6 of the $x_i$ are $\geq 2$ (if 7 or more are $\geq 2$
then $X\geq 2^7=128>71$). If none is, the value is $n\le606$. Say that $x_{i_1}=a_1,\ldots,x_{i_k}=a_k$, $1\le k\le6$, are the variables that are $\ge2$. (Here, $i_1<\cdots<i_k$.)
Define the cyclic distances (gaps) $g_1=i_2-i_1, \ldots , g_{k-1}=i_k-i_{k-1}$, $g_k=i_1-i_k+n$. 
Due to the cyclic nature of $C_{m,5}$, only the values  $a_j\geq 2$ and the gaps $g_j$ are important, and $g_1+\cdots+g_k=n$. 
A monomial $M_i$ contains all variables $\ge2$ exactly when the $4$ variables it leaves out are all equal to $1$. Between $x_{i_j}$ and $x_{i_{j+1}}$ there are $g_j-1$ variables equal to $1$, so there are $\max(g_j-4,0)$ such monomials leaving out $4$ of them.
So we can state that $\sum_{j=1}^k \max(g_j-4,0)$ monomials have value $X=a_1\cdots a_k$.
The squares of the remaining monomials sum to a number $V$ which only depends on $a_1,\ldots, a_k$ and
$\min(g_1,4), \ldots, \min(g_k,4)$. In particular, $V$ does not depend on $m$. 
So summarizing
\[C_{m,5}(x)=V+X^2\sum_{j=1}^k\max(g_j-4,0).\]

So to find all $x$ with $C_{m,5}(x)=9591$ we only need to know $V$ for each possible choice of
\begin{itemize}
\item a tuple $(a_1,\ldots,a_k)$ with $1\le k\le6$, all $a_j\ge2$ and $a_1\cdots a_k\le71$, and
\item a tuple $(\min(g_1,4),\ldots,\min(g_k,4))\in\{1,2,3,4\}^k$.
\end{itemize}
For $k=1,\ldots,6$ there are $70,\,173,\,162,\,65,\,11,\,1$ tuples $(a_1,\ldots,a_k)$, $482$ in total, and each of them comes with $4^k$ tuples $(\min(g_1,4),\ldots,\min(g_k,4))$. This gives
\[70\cdot4+173\cdot4^2+162\cdot4^3+65\cdot4^4+11\cdot4^5+1\cdot4^6=45\,416\]
cases, a quite feasible number. In each of these cases, $C_{m,5}(x)=9591$ means that $9591-V$ is a nonnegative multiple of $X^2$, namely $X^2\sum_{j=1}^k\max(g_j-4,0)$. In fact $9591-V$ is a nonnegative multiple of $X^2$ in only $12$ of the cases (so all other cases are dismissed).
Up to rotation and reflection, these $12$ cases have $k=3$, $(a_1,a_2,a_3)=(2,4,7)$, so $X=56$, and $(g_1,g_2,g_3)=(1,2,3)$ or $(1,3,2)$. So all $g_j\le3$, and then $n=g_1+g_2+g_3=6<16$; equivalently, $V=183$ and $9591-V=3\cdot56^2$ would need $\sum_j\max(g_j-4,0)=3$, while this sum is $0$. Hence $9591\notin\mathcal{C}^*_{m,5}$ for $12\le m\le602$.
\end{proof}

\section{Open problems}\label{s:open}

To make this section self-contained, we repeat some facts from the previous sections. We collect the open problems that seem most interesting or most approachable; all of them are supported by extensive computations.

\begin{prob}\label{pr:145}
Prove C145$'$: $x^2+4y^2+5z^2$ represents every positive integer $4^au$ with $u\equiv6\pmod 8$, except $38$.
\end{prob}
This would make unconditional the universality of the $300$ four-monomial classes (among them \#4 of Section~\ref{4less}) that are now only universal under C145$'$, of rows 49--66 of Table~\ref{tab:settled}, and of $P_{2,4}$ and $C_{2,5}$.

\begin{prob}\label{pr:H}
Is $9591$ the only positive integer not represented by the form $H$ of Section~\ref{s:9591}?
\end{prob}
There is no other miss up to $10^7$. Since every form among the $100$ classes of Section~\ref{s:9591} is a specialization of $H$, $9591$ would then be the only value that is missed by all of them.

\begin{prob}\label{pr:P34}
Are $9591$ and $77871$ the only positive integers not represented by
\[P_{3,4}=(x_1x_2x_3)^2+(x_2x_3x_4)^2+(x_3x_4x_5)^2+(x_4x_5x_6)^2\,?\]
\end{prob}
This holds up to $10^7$ (Section~\ref{s:Ptab}). By Remark~\ref{lem:core} the answer is then the same for all $P_{m,4}$ with $m\ge3$.

\begin{prob}\label{pr:pair}
Are the forms 
\[ x^2+y^2+(xz)^2+(yt)^2+(ztu)^2\]
and 
\[  x^2+y^2+(xz)^2+(yt)^2+(xzu)^2\]
(rows 1 and 8 of Table~\ref{tab:open})
universal?
\end{prob}

Even though these two forms look very similar, the second one is heuristically much more likely to fail universality than the first: the first has $6$ ternary branches and density exponent $e=\tfrac32$, the second has no ternary branch and $e=1$ (Section~\ref{5m5v}).

\begin{prob}\label{pr:row8}
Is the form
\[x^2+y^2+(xz)^2+(yt)^2+(xzu)^2=x^2\bigl(1+z^2+z^2u^2\bigr)+y^2\bigl(1+t^2\bigr)\]
(row 8 of Table~\ref{tab:open}) universal?
\end{prob}
It has no ternary branch, and its number of representations only grows like a power of $\log n$ (Section~\ref{5m5v}). It represents every positive integer up to $10^9$, and up to $10^{10}$ in the residue classes modulo $720$ where it has the fewest representations. A positive answer would also prove that $P_{m,d}$ is universal for all $m\ge3$, $d\ge5$, and $C_{m,d}$ for all $m\ge3$, $d\ge6$ (Proposition~\ref{p:PtableA}).

\begin{prob}\label{pr:P35}
Is $P_{3,5}=(x_1x_2x_3)^2+(x_2x_3x_4)^2+(x_3x_4x_5)^2+(x_4x_5x_6)^2+(x_5x_6x_7)^2$ universal?
\end{prob}
It misses no integer up to $10^7$. It is the smallest form $P_{m,d}$ whose universality is open even assuming C145$'$; a positive answer to Problem~\ref{pr:row8} would settle it.

\begin{prob}\label{pr:C35}
Is $C_{3,5}=(x_1x_2x_3)^2+(x_2x_3x_4)^2+(x_3x_4x_5)^2+(x_4x_5x_6)^2+(x_5x_6x_7)^2+(x_6x_7x_1)^2+(x_7x_1x_2)^2$ universal?
\end{prob}
It misses no integer up to $10^7$. This is a bit surprising, since $C_{m,5}$ misses $9591$ or $77871$ for every $m\ge4$ (Proposition~\ref{p:C5}). Unlike the other open entries of Table~\ref{tab:Cmd}, it does not follow from a positive answer to Problem~\ref{pr:row8}.

\begin{prob}\label{pr:C6}
Are $C_{3,6}$, $C_{4,6}$ and $C_{5,6}$ universal?
\end{prob}
They miss no integer up to $10^7$. Since $C_{m,6}$ has all values of $P_{m,5}$, a positive answer to Problem~\ref{pr:row8} would settle all three, and a positive answer to Problem~\ref{pr:P35} would settle $C_{3,6}$.

\begin{prob}\label{pr:test}
Is there a finite set $S$ of positive integers such that a monomial square form, with any number of monomials, is universal if and only if it represents every element of $S$?
\end{prob}

We have such a set for at most four monomials (Theorem~\ref{t:fifteen}, partly conditional on C145$'$), and conditionally for five monomials in five variables (Proposition~\ref{p:test5}). 
We do not expect such a set to exist when the number of monomials is unbounded, which would make the situation quite different from the 15-theorem for quadratic forms; but a proof of this would be as interesting as finding such a set. 

\clearpage
\appendix

\section{Four monomials}\label{app:4}

Table~\ref{tab:app4} lists the number of classes per number of variables:

\begin{center}\small
\captionof{table}{The $1885$ classes of four-monomial forms by number of variables}\label{tab:app4}
\footnotesize\setlength{\tabcolsep}{3pt}\begin{tabular}{@{}lrrrrrrrrrrrrrr@{}}
\toprule
variables & 3 & 4 & 5 & 6 & 7 & 8 & 9 & 10 & 11 & 12 & 13 & 14 & 15 & total\\
\midrule
all classes & 10 & 66 & 168 & 282 & 363 & 365 & 290 & 187 & 97 & 39 & 13 & 4 & 1 & 1885\\
miss a value $\le143$ & 10 & 62 & 143 & 201 & 197 & 135 & 64 & 22 & 6 & 1 & 0 & 0 & 0 & 841\\
smallest miss $9591$ & 0 & 0 & 1 & 7 & 19 & 28 & 25 & 14 & 5 & 1 & 0 & 0 & 0 & 100\\
univ., unconditional & 0 & 3 & 15 & 44 & 89 & 128 & 136 & 111 & 69 & 32 & 12 & 4 & 1 & 644\\
univ. if C145$'$ & 0 & 1 & 9 & 30 & 58 & 74 & 65 & 40 & 17 & 5 & 1 & 0 & 0 & 300\\
\bottomrule
\end{tabular}
\end{center}

The smallest missed value of the $841$ classes that miss a value $\le143$ (exact computation of all values up to $10^4$):
\begin{center}\small
\captionof{table}{Smallest missed value of the $841$ classes that miss a value $\le143$}\label{tab:app4b}
\begin{tabular}{@{}lrrrrrrrrr@{}}
\toprule
$n$ & 2 & 3 & 6 & 13 & 14 & 15 & 29 & 79 & 143\\
\midrule
classes & 12 & 150 & 56 & 1 & 187 & 235 & 24 & 68 & 108\\
\bottomrule
\end{tabular}
\end{center}

\section{Computations for five monomials}\label{app:5}
This appendix contains the details behind Section~\ref{5m5v}.

\subsection{The $2077$ non-universal classes}

\begin{itemize}
\item \textbf{A small value is missed.} 2056 classes (160\,350 forms) miss some $n\le200$. Smallest missed value:
\end{itemize}
\begin{center}\small
\captionof{table}{Smallest missed value of the $2056$ classes that miss a value $\le200$}\label{tab:app5a}
\begin{tabular}{@{}lrrrrrrrrrrrr@{}}
\toprule
$n$ & 2 & 3 & 6 & 7 & 14 & 15 & 30 & 61 & 78 & 79 & 113 & 143\\
\midrule
classes & 38 & 299 & 336 & 1084 & 4 & 251 & 25 & 2 & 2 & 10 & 2 & 3\\
\bottomrule
\end{tabular}
\end{center}
\begin{itemize}
\item \textbf{Misses above $200$.} Of the $104$ classes that represent every $n\le200$, 10 miss something in $[201,2000]$ ($920$ forms), 3 more in $[2001,6000]$ ($360$), and 6 more in $[6001,10^5]$ ($660$); 2 more ($240$ forms) represent every $n\le10^5$ but miss a larger value. The last two groups are:
\end{itemize}
\begin{center}\small\begin{minipage}{\linewidth}\centering
\captionof{table}{The classes whose smallest missed value is above $6000$}\label{tab:app5b}
\begin{tabular}{rlr}
\toprule
\# & monomials & smallest missed $n$\\
\midrule
1 & $\{x,\,y,\,z,\,xyt,\,xzu\}$ & 7280 \\
2 & $\{x,\,y,\,xz,\,xzt,\,xyzu\}$ & 15839 \\
3 & $\{x,\,y,\,xz,\,xzt,\,yztu\}$ & 15839 \\
4 & $\{x,\,xy,\,xz,\,xyt,\,yzu\}$ & 39759 \\
5 & $\{x,\,xy,\,xz,\,xyt,\,ytu\}$ & 39759 \\
6 & $\{x,\,xy,\,xz,\,yt,\,yztu\}$ & 65663 \\
\midrule
7 & $\{x,\,xy,\,xz,\,yt,\,xytu\}$ & 274647 \\
8 & $\{x,\,xy,\,xz,\,yt,\,xyzu\}$ & 648798 \\
\bottomrule
\end{tabular}\end{minipage}
\end{center}

\subsection{The settled classes}

This subsection lists for each settled class the specializations (certificates) that prove its (conditional) universality (Table~\ref{tab:settled}).
A certificate lists branches $c+Q$ (Definition~\ref{def:branch}); ``$z{=}k$'' fixes a variable and the other variables are free. In rows 1--11, $Q$ is a quaternary form that is universal by the 15-theorem \ref{15}, so $F$ represents every $n\ge c$. In the other rows the exceptional sets of the ternary branches (Section~\ref{Tools}) are compared modulo powers of $2$, using Lemma~\ref{l:cov} and, where possible, the descent of Lemma~\ref{l:scale}; together with an exact search for $n\le60\,000$ this shows that every $n\ge1$ is represented. A certificate with ``subform'' means that the specialization turns $F$ into $c+F'$, where $F'$ is a universal four-monomial form (Theorem~\ref{t:four}). In all cases the finitely many $n<c$ are checked directly.

\small
\begin{longtable}{r p{5.0cm} p{8.3cm}}
\caption{The $66$ settled classes and their certificates}\label{tab:settled}\\
\toprule
\# & $\{M_i\}$ ($F=\sum M_i^2$) & certificate\\
\midrule\endfirsthead
\toprule
\# & $\{M_i\}$ ($F=\sum M_i^2$) & certificate\\
\midrule\endhead
\multicolumn{3}{l}{\emph{Universal by the 15-theorem (11 classes)}}\\
\midrule
1 & $\{x,\,y,\,z,\,t,\,u\}$ & $x{=}0$: $\langle1,1,1,1\rangle$ \\
2 & $\{x,\,y,\,z,\,t,\,xu\}$ & $x{=}1$: $1{+}\langle1,1,1,1\rangle$ \\
3 & $\{x,\,y,\,z,\,t,\,xyu\}$ & $u{=}0$: $\langle1,1,1,1\rangle$ \\
4 & $\{x,\,y,\,z,\,t,\,xyzu\}$ & $u{=}0$: $\langle1,1,1,1\rangle$ \\
5 & $\{x,\,y,\,z,\,t,\,xyztu\}$ & $u{=}0$: $\langle1,1,1,1\rangle$ \\
6 & $\{x,\,y,\,z,\,xt,\,xu\}$ & $x{=}1$: $1{+}\langle1,1,1,1\rangle$ \\
7 & $\{x,\,y,\,z,\,xt,\,tu\}$ & $t{=}1$: $\langle1,1,1,2\rangle$ \\
8 & $\{x,\,y,\,xz,\,xt,\,xu\}$ & $x{=}1$: $1{+}\langle1,1,1,1\rangle$ \\
9 & $\{x,\,y,\,xz,\,zt,\,zu\}$ & $z{=}1$: $\langle1,1,1,2\rangle$ \\
10 & $\{x,\,xy,\,xz,\,xt,\,xu\}$ & $x{=}1$: $1{+}\langle1,1,1,1\rangle$ \\
11 & $\{x,\,xy,\,yz,\,yt,\,yu\}$ & $y{=}1$: $\langle1,1,1,2\rangle$ \\
\midrule
\multicolumn{3}{l}{\emph{Universal, unconditionally (37 classes)}}\\
\midrule
12 & $\{x,\,y,\,z,\,xt,\,yu\}$ & $t{=}0,\,u{=}0$: $\langle1,1,1\rangle$; $t{=}0,\,u{=}1$: $\langle1,1,2\rangle$ \\
13 & $\{x,\,y,\,z,\,xt,\,xyu\}$ & $t{=}0,\,u{=}0$: $\langle1,1,1\rangle$; $t{=}1,\,u{=}0$: $\langle1,1,2\rangle$ \\
14 & $\{x,\,y,\,z,\,xt,\,yzu\}$ & $t{=}0,\,u{=}0$: $\langle1,1,1\rangle$; $t{=}1,\,u{=}0$: $\langle1,1,2\rangle$ \\
15 & $\{x,\,y,\,z,\,xt,\,xtu\}$ & $t{=}0,\,u{=}0$: $\langle1,1,1\rangle$; $t{=}1,\,u{=}0$: $\langle1,1,2\rangle$ \\
16 & $\{x,\,y,\,z,\,xt,\,ytu\}$ & $t{=}0,\,u{=}0$: $\langle1,1,1\rangle$; $t{=}1,\,u{=}0$: $\langle1,1,2\rangle$ \\
17 & $\{x,\,y,\,z,\,xt,\,xyzu\}$ & $t{=}0,\,u{=}0$: $\langle1,1,1\rangle$; $t{=}1,\,u{=}0$: $\langle1,1,2\rangle$ \\
18 & $\{x,\,y,\,z,\,xt,\,xytu\}$ & $t{=}0,\,u{=}0$: $\langle1,1,1\rangle$; $t{=}1,\,u{=}0$: $\langle1,1,2\rangle$ \\
19 & $\{x,\,y,\,z,\,xt,\,yztu\}$ & $t{=}0,\,u{=}0$: $\langle1,1,1\rangle$; $t{=}1,\,u{=}0$: $\langle1,1,2\rangle$ \\
20 & $\{x,\,y,\,z,\,xt,\,xyztu\}$ & $t{=}0,\,u{=}0$: $\langle1,1,1\rangle$; $t{=}1,\,u{=}0$: $\langle1,1,2\rangle$ \\
21 & $\{x,\,y,\,z,\,xyt,\,xyu\}$ & $t{=}0,\,u{=}0$: $\langle1,1,1\rangle$; $x{=}1,\,y{=}1$: $2{+}\langle1,1,1\rangle$ \\
22 & $\{x,\,y,\,xz,\,xt,\,yu\}$ & $x{=}1,\,u{=}0$: $1{+}\langle1,1,1\rangle$; $x{=}1,\,u{=}1$: $1{+}\langle1,1,2\rangle$ \\
23 & $\{x,\,y,\,xz,\,xt,\,zu\}$ & $z{=}1,\,t{=}0$: $\langle1,1,2\rangle$; $x{=}1,\,u{=}0$: $1{+}\langle1,1,1\rangle$ \\
24 & $\{x,\,y,\,xz,\,xt,\,xyu\}$ & $x{=}1,\,u{=}0$: $1{+}\langle1,1,1\rangle$; $x{=}1,\,u{=}1$: $1{+}\langle1,1,2\rangle$ \\
25 & $\{x,\,y,\,xz,\,xt,\,xzu\}$ & $x{=}1,\,u{=}0$: $1{+}\langle1,1,1\rangle$; $x{=}1,\,u{=}1$: $1{+}\langle1,1,2\rangle$ \\
26 & $\{x,\,y,\,xz,\,xt,\,yzu\}$ & $x{=}1,\,u{=}0$: $1{+}\langle1,1,1\rangle$; $x{=}2,\,u{=}0$: $4{+}\langle1,4,4\rangle$ \\
27 & $\{x,\,y,\,xz,\,xt,\,ztu\}$ & $x{=}1,\,u{=}0$: $1{+}\langle1,1,1\rangle$; $x{=}2,\,u{=}0$: $4{+}\langle1,4,4\rangle$ \\
28 & $\{x,\,y,\,xz,\,xt,\,xyzu\}$ & $u{=}0$: $0$ + subform $\{x,\,y,\,xz,\,xt\}$ (universal, Theorem~\ref{t:four}) \\
29 & $\{x,\,y,\,xz,\,xt,\,xztu\}$ & $x{=}1,\,u{=}0$: $1{+}\langle1,1,1\rangle$; $x{=}2,\,u{=}0$: $4{+}\langle1,4,4\rangle$ \\
30 & $\{x,\,y,\,xz,\,xt,\,yztu\}$ & $x{=}1,\,u{=}0$: $1{+}\langle1,1,1\rangle$; $x{=}2,\,u{=}0$: $4{+}\langle1,4,4\rangle$ \\
31 & $\{x,\,y,\,xz,\,xt,\,xyztu\}$ & $u{=}0$: $0$ + subform $\{x,\,y,\,xz,\,xt\}$ (universal, Theorem~\ref{t:four}) \\
32 & $\{x,\,y,\,xz,\,yt,\,zu\}$ & $z{=}1,\,t{=}0$: $\langle1,1,2\rangle$; $z{=}1,\,t{=}1$: $\langle1,2,2\rangle$ \\
33 & $\{x,\,y,\,xz,\,yt,\,xyu\}$ & $y{=}1$: $1$ + subform $\{x,\,t,\,xz,\,xu\}$ (universal, Theorem~\ref{t:four}) \\
34 & $\{x,\,y,\,xz,\,zt,\,tu\}$ & $z{=}0,\,u{=}1$: $\langle1,1,1\rangle$; $z{=}1,\,u{=}0$: $\langle1,1,2\rangle$ \\
35 & $\{x,\,y,\,xz,\,zt,\,xzu\}$ & $z{=}1,\,u{=}0$: $\langle1,1,2\rangle$; $x{=}1,\,z{=}2$: $5{+}\langle1,4,4\rangle$ \\
36 & $\{x,\,y,\,xz,\,zt,\,yzu\}$ & $z{=}1,\,u{=}0$: $\langle1,1,2\rangle$; $z{=}1,\,u{=}1$: $\langle1,2,2\rangle$ \\
37 & $\{x,\,y,\,xz,\,zt,\,xtu\}$ & $z{=}1,\,u{=}0$: $\langle1,1,2\rangle$; $x{=}1,\,t{=}1$: $1{+}\langle1,1,2\rangle$ \\
38 & $\{x,\,y,\,xz,\,zt,\,ztu\}$ & $z{=}1,\,u{=}0$: $\langle1,1,2\rangle$; $z{=}1,\,u{=}1$: $\langle1,2,2\rangle$ \\
39 & $\{x,\,y,\,xz,\,xyt,\,xyu\}$ & $x{=}1$: $1$ + subform $\{y,\,z,\,yt,\,yu\}$ (universal, Theorem~\ref{t:four}) \\
40 & $\{x,\,y,\,xz,\,xzt,\,xzu\}$ & $x{=}1,\,z{=}1$: $2{+}\langle1,1,1\rangle$; $x{=}1,\,z{=}2$: $5{+}\langle1,4,4\rangle$; $x{=}2,\,z{=}1$: $8{+}\langle1,4,4\rangle$ \\
41 & $\{x,\,xy,\,xz,\,xt,\,yu\}$ & $x{=}1,\,u{=}0$: $1{+}\langle1,1,1\rangle$; $x{=}1,\,u{=}1$: $1{+}\langle1,1,2\rangle$ \\
42 & $\{x,\,xy,\,xz,\,xt,\,xyu\}$ & $x{=}1,\,u{=}0$: $1{+}\langle1,1,1\rangle$; $x{=}1,\,u{=}1$: $1{+}\langle1,1,2\rangle$ \\
43 & $\{x,\,xy,\,xz,\,yt,\,yu\}$ & $y{=}1,\,z{=}0$: $\langle1,1,2\rangle$; $x{=}1,\,y{=}2$: $5{+}\langle1,4,4\rangle$ \\
44 & $\{x,\,xy,\,xz,\,yt,\,xyu\}$ & $x{=}1,\,y{=}1$: $2{+}\langle1,1,1\rangle$; $x{=}1,\,y{=}2$: $5{+}\langle1,4,4\rangle$; $x{=}2,\,y{=}1$: $8{+}\langle1,4,4\rangle$ \\
45 & $\{x,\,xy,\,xz,\,xyt,\,xyu\}$ & $x{=}1$: $1$ + subform $\{y,\,z,\,yt,\,yu\}$ (universal, Theorem~\ref{t:four}) \\
46 & $\{x,\,xy,\,yz,\,yt,\,zu\}$ & $y{=}1,\,u{=}0$: $\langle1,1,2\rangle$; $y{=}1,\,u{=}1$: $\langle1,2,2\rangle$ \\
47 & $\{x,\,xy,\,yz,\,yt,\,xyu\}$ & $y{=}1,\,u{=}0$: $\langle1,1,2\rangle$; $x{=}2,\,y{=}1$: $8{+}\langle1,1,4\rangle$ \\
48 & $\{x,\,xy,\,yz,\,yt,\,yzu\}$ & $y{=}1,\,u{=}0$: $\langle1,1,2\rangle$; $y{=}1,\,u{=}1$: $\langle1,2,2\rangle$ \\
\midrule
\multicolumn{3}{l}{\emph{Universal if C145$'$ holds (18 classes)}}\\
\midrule
49 & $\{x,\,y,\,z,\,tu,\,xyt\}$ & $y{=}1,\,t{=}1$: $1{+}\langle1,1,2\rangle$; $y{=}1,\,t{=}2$: $1{+}\langle1,4,5\rangle$ \\
50 & $\{x,\,y,\,z,\,xyt,\,xtu\}$ & $x{=}1,\,t{=}1$: $1{+}\langle1,1,2\rangle$; $x{=}1,\,t{=}2$: $1{+}\langle1,4,5\rangle$ \\
51 & $\{x,\,y,\,xz,\,yt,\,yzu\}$ & $y{=}1,\,z{=}1$: $1{+}\langle1,1,2\rangle$; $y{=}1,\,z{=}2$: $1{+}\langle1,4,5\rangle$ \\
52 & $\{x,\,y,\,xz,\,zt,\,ytu\}$ & $z{=}1,\,u{=}0$: $\langle1,1,2\rangle$; $z{=}2,\,u{=}0$: $\langle1,4,5\rangle$ \\
53 & $\{x,\,y,\,xz,\,zt,\,xyzu\}$ & $z{=}1,\,u{=}0$: $\langle1,1,2\rangle$; $z{=}2,\,u{=}0$: $\langle1,4,5\rangle$ \\
54 & $\{x,\,y,\,xz,\,zt,\,xytu\}$ & $z{=}1,\,u{=}0$: $\langle1,1,2\rangle$; $z{=}2,\,u{=}0$: $\langle1,4,5\rangle$ \\
55 & $\{x,\,y,\,xz,\,zt,\,xztu\}$ & $z{=}1,\,u{=}0$: $\langle1,1,2\rangle$; $z{=}2,\,u{=}0$: $\langle1,4,5\rangle$ \\
56 & $\{x,\,y,\,xz,\,zt,\,yztu\}$ & $z{=}1,\,u{=}0$: $\langle1,1,2\rangle$; $z{=}2,\,u{=}0$: $\langle1,4,5\rangle$ \\
57 & $\{x,\,y,\,xz,\,zt,\,xyztu\}$ & $z{=}1,\,u{=}0$: $\langle1,1,2\rangle$; $z{=}2,\,u{=}0$: $\langle1,4,5\rangle$ \\
58 & $\{x,\,y,\,xz,\,zu,\,xyt\}$ & $z{=}1,\,t{=}0$: $\langle1,1,2\rangle$; $z{=}2,\,t{=}0$: $\langle1,4,5\rangle$ \\
59 & $\{x,\,y,\,xz,\,tu,\,xyt\}$ & $x{=}1,\,t{=}1$: $1{+}\langle1,1,2\rangle$; $x{=}1,\,t{=}2$: $1{+}\langle1,4,5\rangle$ \\
60 & $\{x,\,y,\,xz,\,tu,\,xzt\}$ & $z{=}0,\,u{=}1$: $\langle1,1,1\rangle$; $x{=}1,\,t{=}2$: $1{+}\langle1,4,5\rangle$ \\
61 & $\{x,\,y,\,xz,\,xyt,\,xtu\}$ & $x{=}1,\,t{=}1$: $1{+}\langle1,1,2\rangle$; $x{=}1,\,t{=}2$: $1{+}\langle1,4,5\rangle$ \\
62 & $\{x,\,y,\,xz,\,xzt,\,xtu\}$ & $x{=}1,\,t{=}1$: $1{+}\langle1,1,2\rangle$; $x{=}1,\,t{=}2$: $1{+}\langle1,4,5\rangle$ \\
63 & $\{x,\,xy,\,xz,\,yt,\,tu\}$ & $x{=}1,\,t{=}1$: $1{+}\langle1,1,2\rangle$; $x{=}1,\,t{=}2$: $1{+}\langle1,4,5\rangle$ \\
64 & $\{x,\,xy,\,xz,\,yt,\,xtu\}$ & $x{=}1,\,t{=}1$: $1{+}\langle1,1,2\rangle$; $x{=}1,\,t{=}2$: $1{+}\langle1,4,5\rangle$ \\
65 & $\{x,\,xy,\,xz,\,tu,\,xyt\}$ & $x{=}1,\,t{=}1$: $1{+}\langle1,1,2\rangle$; $x{=}1,\,t{=}2$: $1{+}\langle1,4,5\rangle$ \\
66 & $\{x,\,xy,\,xz,\,xyt,\,xtu\}$ & $x{=}1,\,t{=}1$: $1{+}\langle1,1,2\rangle$; $x{=}1,\,t{=}2$: $1{+}\langle1,4,5\rangle$ \\
\bottomrule
\end{longtable}

\end{document}